\documentclass[12pt]{article}
\usepackage{amssymb}
\usepackage{epsfig}
\usepackage{color}

\normalbaselines

\newtheorem{Theorem}{Theorem}[section]

\newtheorem{Proposition}[Theorem]{Proposition}

\newtheorem{Lemma}[Theorem]{Lemma}

\newtheorem{Corollary}[Theorem]{Corollary}

\newtheorem{Remark}[Theorem]{Remark}

\newtheorem{Example}[Theorem]{Example}

\newcommand{\N}{\mathbb N}

\newcommand{\RR}{{{\rm I} \kern -.15em {\rm R} }}

\newcommand{\C}{{{\rm l} \kern -.42em {\rm C} }}

\newcommand{\nat}{{{\rm I} \kern -.15em {\rm N} }}

\newcommand{\caH}{{\mathcal H}}
\newcommand{\caB}{{\mathcal B}}
\newcommand{\caA}{{\mathcal A}}
\newcommand{\caC}{{\mathcal C}}
\newcommand{\caL}{{\mathcal L}}
\newcommand{\caU}{{\mathcal U}}
\newcommand{\caD}{{\mathcal D}}
\newcommand{\caE}{{\mathcal E}}

\newcommand{\be}{\begin{equation}}
\newcommand{\ee}{\end{equation}}
\newcommand{\beq}{\begin{eqnarray}}
\newcommand{\eeq}{\end{eqnarray}}
\newcommand{\beqs}{\begin{eqnarray*}}
\newcommand{\eeqs}{\end{eqnarray*}}
\newcommand{\bt}{\begin{Theorem}}
\newcommand{\et}{\end{Theorem}}
\newcommand{\br}{\begin{Remark}}
\newcommand{\er}{\end{Remark}}
\newcommand{\bc}{\begin{Corollary}}
\newcommand{\ec}{\end{Corollary}}
\newcommand{\bl}{\begin{Lemma}}
\newcommand{\el}{\end{Lemma}}
\newcommand{\bd}{\begin{definition}}
\newcommand{\ed}{\end{definition}}

\title{Exponential stability of abstract evolution equations with
 time delay}
\author{{\sc Serge Nicaise}
\\Universit\'e de Valenciennes et du Hainaut Cambr\'{e}sis\\
ISTV--LAMAV\\
59313 Valenciennes Cedex 9, France
\\\\
{\sc Cristina Pignotti}
\\Dipartimento di Ingegneria e Scienze dell'Informazione e Matematica\\
 Universit\`{a} di L'Aquila\\
Via Vetoio, Loc. Coppito, 67010 L'Aquila Italy}

\date{}

\begin{document}

\textwidth=160 mm

\textheight=225mm

\parindent=8mm

\frenchspacing

\maketitle

\begin{abstract}
We consider abstract semilinear evolution equations with a time delay feedback.
We show that, if the $C_0$-semigroup describing the linear part of the model is exponentially stable, then the whole system retains this good property when a suitable {\it smallness} condition on the time delay feedback is satisfied.
Some  examples illustrating our abstract approach are also given.
\end{abstract}

\vspace{5 mm}

\def\qed{\hbox{\hskip 6pt\vrule width6pt
height7pt
depth1pt  \hskip1pt}\bigskip}

%% {\bf 2000 Mathematics Subject Classification:}
%%35L05, 93D15

 %%{\bf Keywords and Phrases:}  wave equation,  delay feedbacks, stabilization

\section{Introduction}
\label{intro}\hspace{5mm}

\setcounter{equation}{0}

Let $\caH$ be a fixed Hilbert space with norm $\|\cdot\|$, and consider an operator $\caA$ from $\caH$ into itself that generates a $C_0$-semigroup  $(S(t))_{t\geq 0}$ that is exponentially stable, i.e.,
there exist two positive constants $M$ and $\omega$ such that
\begin{equation}\label{assumpexpdecay}
\|S(t)\|_{{\mathcal L}(\caH)}\leq M e^{-\omega t},\  \forall t\geq 0,
\end{equation}
where, as usual, ${\mathcal L}(\caH)$ denotes the space of  bounded linear 
operators from $\caH$ into itself.
For a fixed delay parameter $\tau$, a fixed bounded operator $\caB$ from $\caH$ into itself
and for a real parameter $k$,
we consider the evolution equation

 \begin{eqnarray}\label{abstract}
\left\{
\begin{array}{ll}
& U_{t}(t) = \caA U(t)+F(U(t))+k \caB  U(t-\tau) \quad \mbox{\rm in }
(0,+\infty)\\
& U(0)=U_0,\  \caB U(t-\tau)=f(t), \quad \forall t\in (0,\tau),
\end{array}
\right.
\end{eqnarray}
where $F: \caH\rightarrow\caH$  satisfies some Lipschitz conditions, the initial datum $U_0$ belongs to $\caH$ and $f\in C([0,\tau]; \caH)$.

Time delay effects often appear in many applications and physical problems. On the other hand, it is well-known (cfr. \cite{Batkai,Datko,DLP,NPSicon06,XYL})
that they can induce some instability.
Hence we are interested in giving an exponential stability result for such a problem under a suitable condition between
the constant  $k$ and the constants $M, \omega, \tau,$ the norm of $\caB$ and the nonlinear term $F.$
For some particular examples (see e.g. \cite{ANP10, Batkai, SCL12, guesmia, Said, AlNP})
we know that the above problem, under certain {\sl smallness} conditions on the delay feedback $k\caB,$
is exponentially stable,
the proof being from time to time quite technical
because some observability inequalities or perturbation methods are used.
Hence our main goal is to furnish a direct proof of this stability result by using the so-called
  Duhamel's formula (or variation of parameters formula).

Observe that our proof is simpler with respect to the ones used so far for particular models. Moreover, we emphasize its generality. Indeed, it applies to every model in the form
(\ref{abstract}) when the operator $\caA$ generates an exponentially stable semigroup.

In the same spirit, we want to prove existence and exponential stability results
when the operator $\caB$ is unbounded (and $F=0$). In that case  they are proved using  Duhamel's formula
but under some admissibility conditions.

Note also that previous papers deal with linear models and $\caB$  bounded, while here  we include a nonlinear term $F$ or $\caB$  unbounded.

The paper is organized as follows. In section
2 we study the case with bounded feedback operator $\caB$
and nonlinear term $F$ Lipschitz, giving a well-posedness result and an exponential decay estimate. The analysis is then extended to a more general linear 
term $F$ in section 3, under more restrictive assumptions. In section
4, we consider, only for the linear model, unbounded delay feedback
operators $\caB$ and prove a well-posedness and an exponential stability result. Finally, in section 5, some  illustrative examples with $\caB$ unbounded are given.

\section{The case $F$ globally Lipschitz}

\label{pbform}\hspace{5mm}

\setcounter{equation}{0}

In this section, we assume that $F$ is globally Lipschitz continuous, namely
\begin{equation}\label{lip}
\exists \gamma >0 \quad \mbox{\rm such\ that}\quad \Vert F(U_1)-F(U_2)\Vert_\caH
\le \gamma \Vert U_1-U_2\Vert_\caH,\quad \forall\ U_1, U_2\in \caH\,.
\end{equation}
Moreover, we assume that $F(0)=0.$

The following well--posedness result holds.

\begin{Proposition}\label{WellP}
For any initial datum $U_0\in {\mathcal H}$ and $f\in C ([0,\tau]; \caH)$, there exists a unique (mild) solution
$U\in C([0,+\infty), {\mathcal H})$  of problem $(\ref{abstract}).$
Moreover,
\begin{equation}\label{Du}
U(t)=S(t)U_0+\int_0^t S(t-s)[F(U(s))+ k\caB U(s-\tau )]\,ds.
\end{equation}
\end{Proposition}

\noindent{\bf Proof.}
We use an iterative argument.
Namely in the interval $(0,\tau)$, problem (\ref{abstract}) can be seen as an inhomogeneous evolution problem
\begin{eqnarray}\label{abstract2}
\left\{
\begin{array}{ll}
& U_{t}(t) = \caA U(t)+F(U(t))+g_0(t)\quad \mbox{\rm in }
(0,\tau)\\
& U(0)=U_0,
\end{array}
\right.
\end{eqnarray}
where $g_0(t)=k f(t)$. This problem has a unique
solution
$U\in  C([0,\tau], {\mathcal H})$ (see Th. 1.2, Ch. 6 of
\cite{pazy})
satisfying
$$
U(t)=S(t)U_0+\int_0^t S(t-s)[F(U(s))+ g_0(s)]\,ds.
$$
This yields $U(t),$ for $t\in [0,\tau].$
Therefore on $(\tau,2\tau)$, problem (\ref{abstract}) can be seen as an inhomogeneous evolution problem

\begin{eqnarray}\label{abstract3}
\left\{
\begin{array}{ll}
& U_{t}(t) = \caA U(t)+F(U(t))+g_1(t)\quad \mbox{\rm in }
(\tau, 2\tau)\\
& U(\tau)=U(\tau-),
\end{array}
\right.
\end{eqnarray}
where $g_1(t)=k \caB  U(t-\tau)$. Hence, this problem has a unique solution
$U\in C([\tau,2\tau], {\mathcal H})$ given by
$$
U(t)=S(t-\tau)U(\tau-)+\int_\tau^t S(t-s) [F(U(s))+g_1(s)]\,ds, \forall t\in [\tau,2\tau].
$$
By iteration, we obtain a global solution $U$ satisfying (\ref{Du}).\qed

Now we will prove the following exponential stability result.

\begin{Theorem}\label{stab2}
Let $M, \omega, \gamma$ as in
 $(\ref{assumpexpdecay})$ and $(\ref{lip}).$
There is a positive constant $k_0$ such that for  $k$ satisfying
\begin{equation}\label{condkpetit}
\vert k\vert <k_0:=\frac{e^{\tau\omega }-1}{\tau \|\caB\|_{{\mathcal L}(\caH)} M e^{\tau \omega }},
\end{equation}
and for $\gamma < \gamma (\vert k\vert ),$
where $\gamma (\vert k\vert )$ is a suitable constant depending on $\vert k\vert,$ then
 there exist $\omega' >0$ and $M'>0$ such that
the solution $U\in C([0,+\infty), {\mathcal H})$  of problem $(\ref{abstract}),$
with $U_0\in {\mathcal H}$ and $f\in C([0,\tau ]; \caH),$ satisfies
\begin{equation}\label{exponentiald}
\|U(t)\|_\caH\le M' e^{-\omega' t} (\|U_0\|_\caH+\int_0^\tau   e^{\omega s}  \|f(s)\|_\caH\,ds),
\quad \forall t\ge \tau.
\end{equation}
From its definition the constant $k_0$ depends only on $M,  \omega, \tau$ and the norm of $\caB$.
\end{Theorem}
\noindent{\bf Proof.}
First assume $F\equiv 0.$
We use again an iterative argument and  Duhamel's formula but here on the whole ${\mathbb R}_+$, namely
we can write
\begin{equation}\label{Duhamelformula}
U(t)=S(t)U_0+k\int_0^t S(t-s)  \caB  U(s-\tau)\,ds,\ \  \forall t>0.
\end{equation}

Then,
\[
\Vert U(t)\Vert_\caH \leq M e^{-\omega t} (\|U_0\|_\caH+|k|\int_0^t   e^{\omega s}  \|\caB  U(s-\tau)\|_\caH\,ds),\ \  \forall t>0.
\]

Let us show that this implies that
\begin{equation}\label{estiterativestar}
e^{\omega t}
\|U(t)\|_\caH\leq  M (\|U_0\|_\caH+|k|\alpha)
 (1+|k| \tau B M e^{\omega \tau})^n, \quad t \in [0,(n+1)\tau],\ n\in \N,
\end{equation}
where for shortness  we have set $\alpha:=\int_0^\tau   e^{\omega s}  \|f(s)\|_\caH\,ds$ and $B=\|\caB\|_{{\mathcal L}(\caH)}$.

Now  if  we set
\begin{equation}\label{defsigma}
\sigma=\tau^{-1} \ln(1+|k| \tau B M e^{\omega \tau}),
\end{equation}
then we see that  (\ref{estiterativestar}) gives
$$
\| U(t)\|_\caH\leq M e^{-(\omega -\sigma) t} (\|U_0\|_\caH+|k|\alpha), \forall t>0.
$$
Hence $\|U(t)\|_\caH$ will decay exponentially if
$\sigma -\omega$ is negative or equivalently if
\[
1+|k| \tau B M e^{\omega \tau}<e^{\tau \omega},
\]
which is nothing else than (\ref{condkpetit}).

Under  this constraint,
we deduce that the  estimate  (\ref{exponentiald}) holds with $\omega'=
\omega -\sigma.$

Hence we are reduced to prove (\ref{estiterativestar}).
First in $(0,\tau)$, (\ref{Du}) and the initial condition from (\ref{abstract}) yield
\[
\|U(t)\|_\caH\leq  M e^{-\omega t} (\|U_0\|_\caH+|k|\int_0^t   e^{\omega s}  \|f(s)\|_\caH\,ds),\ \  \forall t\in (0,\tau).
\]

Then,
\[
e^{\omega t}\|U(t)\|_\caH\leq  M  (\|U_0\|_\caH+|k|\alpha ),\ \  \forall t\in (0,\tau),
\]
which is nothing else than (\ref{estiterativestar}) for $n=0$.

Second for any $m\in \N^*$, we assume that (\ref{estiterativestar}) holds for all $n\leq m-1$
and prove that it holds for $m$. Indeed by  (\ref{Du}), we have for all $t\in (m\tau, (m+1)\tau)$
$$
\|U(t)\|_\caH\leq M e^{-\omega t} \Big(\|U_0\|_\caH+|k| \alpha
+|k|
\sum_{\ell=0}^{m-1}\int_{(\ell+1)\tau}^{(\ell+2)\tau}   e^{\omega s}  \|\caB  U(s-\tau)\|_\caH\,ds\Big).
$$
Now for $s\in ((\ell+1)\tau, (\ell+2)\tau)$, with $\ell=0, \ldots, m-1$, we notice that
$s-\tau$ belongs to $ (\ell\tau, (\ell+1)\tau)$, and  using our iterative assumption, we get

\begin{eqnarray*}
\|U(t)\|_\caH&\leq& M e^{-\omega t} \Big(\|U_0\|_\caH+|k| \alpha
\\
&+&|k|
\sum_{\ell=0}^{m-1}\int_{(\ell+1)\tau}^{(\ell+2)\tau}   e^{\omega s}
B M e^{-\omega (s-\tau)} (\|U_0\|_\caH+|k|\alpha) (1+|k| \tau B M
e^{\omega\tau})^\ell
\,ds\Big)
\\
&\leq& M e^{-\omega t} \Big(\|U_0\|_\caH+|k| \alpha
+|k| \tau B M (\|U_0\|_\caH+|k|\alpha) e^{\omega \tau}
\sum_{\ell=0}^{m-1}
 (1+|k| \tau B M e^{\omega \tau})^\ell\Big)
\end{eqnarray*}
Hence  we have obtained that

$$
\|U(t)\|_\caH \leq  M e^{-\omega t} (\|U_0\|_\caH+|k|\alpha)
\Big(1+|k| \tau B M  e^{\omega \tau}
\sum_{\ell=0}^{m-1}
 (1+|k| \tau B M e^{\omega \tau})^\ell\Big).
$$
Because one readily checks that
\[
1+|k| \tau B M  e^{\omega \tau}
\sum_{\ell=0}^{m-1}
 (1+|k| \tau B M e^{\omega \tau})^\ell
=
 (1+|k| \tau B M e^{\omega \tau})^m,
\]
we obtain
$$
e^{\omega t}
\|U(t)\|_\caH\leq  M (\|U_0\|_\caH+|k|\alpha)
 (1+|k| \tau B M e^{\omega \tau})^m, \quad t\in \ [m\tau, (m+1)\tau ].
$$
This estimate and the recurrence assumption,  as
$$(1+|k| \tau B M e^{\omega \tau})^n\leq
(1+|k| \tau B M e^{\omega \tau})^m
 \quad\mbox{\rm for\ all}
 \ \  n\leq m-1,$$
 imply that
(\ref{estiterativestar}) holds for $m$. 
So the result is proved in the linear case. 

In order to extend it to the nonlinear model,
let us introduce (cfr. \cite{NPSicon06}) the new variable

$$Z(t,\rho ):=\caB U(t-\tau\rho ),\quad \rho\in (0,1),\ t>0\;.$$
Then problem (\ref{abstract}) may be rewritten as
\begin{equation}\label{riscritto}
\left\{
\begin{array}{l}
U_{t}(t) = \caA U(t)+F(U(t))+k Z(t,1)\quad \mbox{\rm in }
(0,+\infty)\\
Z_t(t, \rho )=-\tau^{-1}Z_{\rho}(t,\rho )\quad
 \mbox{\rm in }
(0,+\infty )\times (0,1)\\
Z(t,0)=\caB U(t)\\
 U(0)=U_0,\  Z(0, \rho )=f((1-\rho) \tau ), \quad \forall \rho\in (0,1).
\end{array}
\right.
\end{equation}
Therefore, if we set $V:=(U, Z)^T,$ the linear part of (\ref{riscritto}), namely

$$
\left\{
\begin{array}{l}
U_{t}(t) = \caA U(t)+k Z(t,1)\quad \mbox{\rm in }
(0,+\infty)\\
Z_t(t, \rho )=-\tau^{-1}Z_{\rho}(t,\rho )\quad
 \mbox{\rm in }
(0,+\infty )\times (0,1)\\
Z(t,0)=\caB U(t)\\
 U(0)=U_0,\  Z(0, \rho )=f((1-\rho ) \tau ), \quad \forall \rho\in (0,1),
\end{array}
\right.
$$
becomes
$$
\left\{
\begin{array}{l}
V_t=\tilde {\caA} V\\
V(0)=(U(0), Z(0, \cdot ))^T.
\end{array}\right.
$$
It is easy to see that $\tilde {\caA} $ generates a strongly continuous semigroup $(T(t))_{t\ge 0}$ in
the Hilbert space
$\tilde \caH := \caH\times L^2(0,1;\caH ).$
Moreover, $(T(t))_{t\ge 0}$ is exponentially stable.

Indeed, we clearly have
$$
\Vert Z(t,\rho) \Vert^2_{L^2(0,1;\caH )}=\int_0^1 \Vert \caB U(t-\tau\rho )\Vert_{\caH}^2d\rho\le B^2\int_0^1\Vert U(t-\tau\rho )\Vert_{\caH}^2d\rho\,.
$$
Then, the exponential estimate for $U$ gives, for $t\ge 2\tau,$

$$
\Vert Z(t,\rho) \Vert^2_{L^2(0,1;\caH )}\le  B^2{M^\prime}^2h^2e^{-2\omega^\prime t}
\Big (   \Vert U_0\Vert_{\caH}  +\int_0^\tau e^{\omega s} \Vert f(s)\Vert_{\caH} ds \Big )^2
$$
with $h^2:=\frac 1 {\tau }\int_0^\tau e^{2\omega^\prime s} ds.$
Thus, there exists a positive constant $\tilde M,$ depending on $M, \omega, \tau, \vert k\vert $ and the norm of $\caB,$ such that
\begin{equation}\label{stimasecondsemig}
\Vert T(t)\Vert_{\tilde \caH }\le \tilde M e^{-\omega^\prime  t},\quad t>0.
\end{equation}
Coming back to (\ref{riscritto})
and using Duhamel's formula,
$V:=(U, Z)^T$ can be written as
$$
V(t) =T(t) V_0+ \int_0^t T(t-s) \tilde F (V(s)) ds,$$
where $\tilde F (V(s))= (F(U(s)),0)^T.$
Therefore,
$$\Vert V(t)\Vert_{\tilde \caH}
\le \tilde Me^{-\omega^\prime t} \Vert V_0\Vert_{\tilde\caH }+
 \tilde Me^{-\omega^\prime t} \int_0^t e^{\omega s }\gamma
\Vert V(s)\Vert_{\tilde \caH}ds 
$$
and the exponential stability estimate follows from Gronwall's lemma if
$\omega^\prime -\gamma \tilde M <0.$
$\qed $

\begin{Remark}
{\rm
From our proof we see that for $F\equiv 0$ the explicit decay of $\|U(t)\|_\caH$
is
\begin{equation}\label{nu1}
\|U(t)\|_\caH\leq C  e^{(\sigma-\omega) t} (\|U_0\|_\caH+\alpha),
\end{equation}
for some $C>0$.
}
\end{Remark}

\begin{Remark} {\rm
 Note that Theorem \ref{stab2} is very general. Indeed, it gives stability results, when the delay feedback parameter $k$ is sufficiently small, for every model in the form (\ref{abstract}) if the semigroup $(S (t))_{t\ge 0}$ generated by the linear operator ${\caA }$ is exponentially stable.
For instance, it furnishes stability results for previously studied models for wave type equations (cfr. \cite{ANP10, SCL12, MCSS}), Timoshenko models (cfr.
\cite{Said}). Also, it includes recent stability results for problems with
viscoelastic damping and time delay (cfr. \cite{guesmia, AlNP}).}
\end{Remark}

\section{More general nonlinearities}

\subsection{Abstract existence and stability results}
Here we consider a more general class of nonlinearities. More precisely, we assume that
for every constant $c$ there exists a  positive constant $L(c)$ such that
\begin{equation}\label{lip1}
\Vert F(U_1)-F(U_2)\Vert_\caH \le L(c)\Vert U_1-U_2\Vert_\caH\,,
\end{equation}
for all $U_1, U_2\in \caH$ with $\Vert U_1\Vert_\caH \le c,\  \Vert U_2\Vert_\caH \le c.$

Moreover, we assume that there exists an increasing continuous function
$\chi:[0,+\infty)\rightarrow [0,+\infty),$ with $\chi (0)=0,$ such that
\begin{equation}\label{lip2}
\Vert F(U)\Vert_\caH \le\chi (\Vert U\Vert_\caH)\Vert U\Vert_\caH, \ \forall\ U\in \caH\,.
\end{equation}

Now, the nonlinear term introduces additional difficulties. We can give an exponential stability result under a well--posedness assumption for {\sl small} initial data.
Then, we will show that this assumption is satisfied for a quite large class of examples.

\begin{Theorem}\label{Wellnonlinear}
Let $\tilde M, \omega^\prime$ be as in $(\ref{stimasecondsemig}).$
Suppose that for $\vert k\vert$ sufficiently small
\begin{equation}\label{assump}
\begin{array}{lll}
\exists\ \rho_0>0 \ \mbox{\rm   and  } \  C_{\rho_0}>0\ \mbox{\rm such  that}\ \forall
\ U_0\in \caH ,\ f\in C([0,\tau];\caH )\\\medskip
 \mbox{\rm with}\
(\Vert U_0\Vert_\caH^2 +\int_0^\tau \vert k\vert \Vert f(s)\Vert_\caH^2 ds )^{1/2}<\rho_0,
\
\mbox{\rm there\ exists\ a\ unique}\ \mbox{\rm global\ solution} \\\medskip U \in C([0,+\infty, \caH) \
\mbox{\rm to}\ (\ref{abstract})
\ \mbox{\rm with}\ \ \Vert U(t)\Vert_\caH\le C_{\rho_0}<\chi^{-1}\left(\frac {\omega^\prime}{\tilde M }\right ),\ \forall\ t>0.
\end{array}
\end{equation}
Then there exists $\tilde k>0$ such that if $\vert k\vert <\tilde k,$ for every $U_0\in\caH$ and $f\in C([0,\tau];\caH )$ satisfying the assumption from $(\ref{assump}),$ the solution $U$ of problem $(\ref{abstract})$ satisfies the exponential decay estimate
 \begin{equation}\label{exponentialdbis}
\|U(t)\|_\caH\le  M^* e^{-\tilde \omega t} (\|U_0\|_\caH+\int_0^\tau   e^{\omega s}  \|f(s)\|_\caH\,ds),
\quad \forall t\ge \tau,
\end{equation}
for suitable constants $M^*,\ \tilde\omega.$
 \end{Theorem}

\noindent {\bf Proof.}  We can simply repeat the  previous  proof of Theorem \ref{stab2} with $\chi(C_{\rho_0})$
instead of $\gamma.$\qed

\subsection{Examples}\label{esempi}

We now give some examples for wich assumption (\ref{assump}) is satisfied.

Let $H$ be a real Hilbert space,  with norm $\|\cdot\|_H,$ and let
$A_1:\caD(A_1)\rightarrow H,$ a positive self--adjoint operator with a compact inverse in $H.$ Denote by $V:=\caD(A_1^{\frac 12})$ the domain of $A_1^{\frac 12}.$ Further, for i=1,2, let $W_i$ be a real Hilbert space (which will be identified to its dual space) and let $C\in {\mathcal L}(W_1,H),$ $B\in {\mathcal L}(W_2, H).$
Assume that, for some constant $\mu >0$
\begin{equation}\label{contrast}
\Vert B^*u\Vert^2_{W_2}\le\mu \Vert C^*u\Vert^2_{W_1},\quad \forall\ u\in V\,.
\end{equation}
Let be given a functional $G:V\to \RR$ such that $G$ is G\^ateaux differentiable at any $x\in V$. We further assume (cfr. \cite{ACS}) that
\\a) For any $u\in V$   there exists a positive constant $c(u)$ such that
\[
|DG(u)(v)|\leq c(u) \|v\|_H,\  \forall v\in V,
\]
where $DG(u)$ is the G\^ateaux derivative of $G$ at $u$. Consequently $DG(u)$ can be extended in the whole $H$ and we will denote by $\nabla G(u)$ the unique element in $H$
such that
\[
(\nabla G(u),v)_H=DG(u)(v),\  \forall v\in H.
\]
b) For all $c>0$, there exists $L(c)>0$ such that
\[
\|\nabla G(u)-\nabla G(v)\|_H\leq L(c)\|A_1^{\frac12} (u-v)\|_H
\]
for all $u, v\in V$ such that $\|A_1^{\frac12} u\|_H\leq c$ and $\|A_1^{\frac12} v\|_H\leq c$.
\\
c) There exists a suitable increasing continuous function $\psi$  satisfying $\psi (0)=0$ such that

$$\|\nabla G(u)\|_H\le \psi (\Vert A_1^{\frac 12} u\Vert )\Vert A_1^{\frac 12} u\Vert_H^2,
\  \forall u\in V.$$

 In this setting let us consider the second order evolution equation
\begin{equation}\label{second}
\begin{array}{l}
u_{tt}+A_1u+CC^*u_t=\nabla G(u)+kBB^*u_t(t-\tau),\quad t>0, \\
u(0)=u_0,\ u_t(0)=u_1, \ B^*u_t(t-\tau)=g(t),\ t\in (0,\tau),
\end{array}
\end{equation}
with $(u_0,u_1)\in V\times H$.
Denoting $v:=u_t$ and $U:= (u, v)^T,$ this problem may be rewritten in the form
(\ref{abstract}) with
$$\caA:=\left (
\begin{array}{l}
\ \ 0\quad \ \quad\ 1\\
-A_1\ \ -CC^*
\end{array}
\right ),
$$
$$F(U):=(0, \nabla G(u))^T,\quad \caB U:= (0, BB^*v)^T\,.$$

\noindent
The above assumptions on $G$ imply that $F$ satisfies (\ref{lip1}) and (\ref{lip2})  in $\caH := V\times H,$
with $\chi=\psi.$
We define the energy of solutions of problem (\ref{second}) as

\begin{equation}\label{energy}
E(t):=E(t,u(\cdot))=\frac 1 2 \Vert u_t\Vert_H^2+\frac 1 2\Vert A_1^{\frac 12}
u\Vert_H^2- G(u)+\frac 12 \int_{t-\tau}^t\vert k\vert \Vert B^*u_t(s)\Vert_{W_2}^2 ds\,.
\end{equation}

\noindent
We will show that for the above model Theorem \ref{Wellnonlinear} holds.

First of all note that
$$
\begin{array}{l}
\displaystyle{
E^\prime (t)=-\Vert C^*u_t(t)\Vert^2_{W_1}+k\langle B^*u_t(t),B^*u_t(t-\tau)\rangle
+\frac {\vert k\vert} 2\Vert B^*u_t(t)\Vert^2_{W_2}-\frac{\vert k\vert} 2
\Vert B^* u_t(t-\tau)\Vert^2_{W_2}}\\
\hspace{1 cm}\le -\Vert C^* u_t(t)\Vert_{W_1}^2+ \vert k\vert\,
\Vert B^*u_t(t)\Vert^2_{W_2}
\end{array}
$$
Then, if $\vert k\vert <\frac 1{\mu},$ the energy is not increasing.

We can prove the following well-posedness result for sufficiently small data.

\begin{Proposition}
The assumption $(\ref{assump})$ is satisfied for $\vert k\vert<1/\mu$.
\end{Proposition}

\noindent {\bf Proof.}
Note that the condition $\vert k\vert <1/\mu$ guarantees that  the energy is not increasing.

First of all, on $[0,\tau ]$ the abstract system   may be rewritten in the form (\ref{abstract2}) with
$g_0(t)=(0, kBg(t)).$

Then, from classical theory for nonlinear evolution equation (see Th. 1.4, Ch. 6 of \cite{pazy}), there exists a unique mild local solution $U$ defined in
a maximal time interval  $[0, \sigma)$ with $0<\sigma\le\tau.$
We will show that $\sigma =\tau\,.$

We argue similarly to \cite{ACS}.
Note that if $\psi (\Vert A_1^{\frac 12}u_0\Vert_H )<\frac 1 4,$ then
$$E(0)\ge \frac 1 2 \Vert u_1\Vert_H^2+\frac 1 2  \Vert A_1^{\frac 12}u_0\Vert_H^2
-G(u_0)\ge \frac 1 2 \Vert u_1\Vert_H^2+\frac 1 4  \Vert A_1^{\frac 12}u_0\Vert_H^2\,\ge 0\,.$$
We first show that if
\begin{equation}\label{condA1}
\psi (\Vert A_1^{\frac 12}u_0\Vert_H )<\frac 1 4 \quad \hbox{
and } \quad \psi (2(E(0))^{\frac 12})<\frac 14,
\end{equation}
then
\begin{equation}\label{A1}
E(t)\ge \frac 1 2 \Vert u_t(t)\Vert_H^2+\frac 1 4  \Vert A_1^{\frac 12}u(t)
\Vert_H^2\,,\  \forall t\in [0,\sigma ).
\end{equation}

Let $r:=\sup \{s\in [0,\sigma) $ such that (\ref{A1}) holds for every $t\in [0,s]\}.$ Suppose that $r<\sigma,$ then,
\begin{equation}\label{A2}
E(r)\ge \frac 1 2 \Vert u_t(r)\Vert_H^2+\frac 1 4  \Vert A_1^{\frac 12}u(r)
\Vert_H^2\ge 0\,.
\end{equation}

Thus, from (\ref{A2}), we have
$$\psi (\Vert A_1^{\frac 12}u(r)\Vert_H )\le \psi (2(E(r))^{\frac 12}) <
\psi  (2(E(0))^{\frac 12}) <\frac 1 4\,.
$$
This gives
$$E(r)\ge  \frac 1 2 \Vert u_t(r)\Vert_H^2+\frac 1 2  \Vert A_1^{\frac 12}u(r)
\Vert_H^2-G(u(r)) >  \frac 1 2 \Vert u_t(r)\Vert_H^2+\frac 1 4  \Vert A_1^{\frac 12}u(r)
\Vert_H^2,$$
which contradicts the maximality of $r.$
This implies $r=\sigma.$

Now, let us set
$$\rho_0= \min \left\{\frac 12 \psi^{-1}\left (\frac 1 4\right ), \frac 1 {2\sqrt{2}} \psi^{-1}\left (  \frac {\omega^\prime} {\tilde M}\right )\right\}>0.$$
In a second step we show that (\ref{condA1}) holds   for all $u_0\in \caD (A_1^{\frac 12})$, $u_1\in H,$ $g\in C([0,\tau], W_2),$
satisfying
\begin{equation}\label{condexistenceglobale}
\Big (\Vert A_1^{\frac 12} u_0\Vert_H^2+\Vert u_1\Vert_H^2+
\int_0^{\tau}\vert k\vert \Vert g(s)\Vert_{W_2}^2 ds
\Big )^{\frac 12}<\rho_0\,.\end{equation}
Indeed, as this assumption implies that  $\Vert A_1^{\frac 12} u_0\Vert_H <\rho_0$, then one has
$$\psi (\Vert A_1^{\frac 12} u_0\Vert_H )<\psi (\rho_0)=\psi (\frac 12\psi^{-1}(\frac 14))< \frac 14.$$
Hence by the assumption c) on $G$, we deduce that
$$E(0)\le \frac 34\Vert A_1^{\frac 12} u_0\Vert_H^2+\frac 12\Vert u_1\Vert_H^2+
\frac 12 \int_0^{\tau}\vert k\vert \Vert g(s)\Vert_{W_2}^2 ds< \rho_0^2,$$
and, by definition of $\rho_0,$ we conclude that
$$\psi(2 (E(0))^{\frac 12})<\psi (\psi^{-1}(\frac 14))=\frac 14\,.$$

In conclusion under the assumption (\ref{condexistenceglobale}), the estimate (\ref{condA1}) holds,
implying in particular that
$$
0\leq \frac 1 2 \Vert u_t(t)\Vert_H^2+\frac 1 4  \Vert A_1^{\frac 12}u(t)
\Vert_H^2\leq E(t)\le E(0)\le
\rho_0^2, \forall t\in [0,\sigma].$$
 Then again by   \cite[Th. 1.4, Ch. 6]{pazy}), $\sigma=\tau.$

Now we can consider the interval $[\tau, 2\tau )$ and we can rewrite the problem in the form (\ref{abstract2}) with $g_1(s)=(0, kBu_t(t-\tau)).$
As before, there exists a local solution and   arguing as on $[0,\tau ]$ we  obtain  a solution
on $[0, 2\tau ]$ under the assumption  (\ref{condexistenceglobale}).

By repeating this argument we prove that, if
(\ref{condexistenceglobale}) holds,
 then the solution exists on $[0, +\infty)$ and
$$\Vert u_t(t)\Vert_H^2+\Vert A_1^{\frac 12}u(t)\Vert_H^2<4\rho^2_0\le \frac 1 2 \left [
\psi^{-1}\left (\frac {\omega^\prime} {\tilde M}  \right )
\right ]^2
.$$
This proves (\ref{assump}).\qed

If $\caA$ generates an exponentially stable continuous semigroup on $\caH,$ then the exponential estimate (\ref{exponentialdbis}) holds for $k$ small enough,
for small initial data.

\begin{Remark}{\rm
The abstract model (\ref{second}) includes semilinear versions
of previously analyzed concrete models for wave-type equations
(cfr. \cite{NPSicon06}); see the example below. Of course, due to the presence of the nonlinearity
we obtain the stability result (for small initial data) under a more restrictive assumption on the size
of the delay feedback parameter $k$.
Observe also that  models with viscoelastic damping could be considered but
with also an extra not--delayed damping necessary to avoid blow-up of
solutions, at least for small data.}
\end{Remark}

\begin{Example}
{\rm
As an explicit example of system (\ref{second}) let us consider the wave equation with local internal damping
and internal delay.
More precisely, let $\Omega\subset\RR^n,$ $n\ge 3,$ be an open bounded domain   with a boundary
$\partial\Omega$ of class $C^2.$  
Denoting by $m$ the standard multiplier 
$m(x)=x-x_0,\ x_0\in\RR^n,$ 
let $\omega_1$ be the intersection of $\Omega$ with an open neighborhood of the
subset of $\partial\Omega$
\be\label{defgamma0}\Gamma_0=\{\, x\in\partial\Omega\, :\ m(x)\cdot \nu (x)>0\,\}.
\ee
Moreover, let $\omega_2$ be any set satisfying $\omega_2\subseteq\omega_1.$
Let us consider the initial boundary value  problem
\begin{eqnarray}
& &u_{tt}(x,t) -\Delta u (x,t)+a\chi_{\omega_1} u_t(x,t)+k\chi_{\omega_2}
 u_t(x,t-\tau)\nonumber\\\medskip
& & \hspace{6.5 cm}=
\vert u(x,t)\vert^{\beta} u(x,t)\quad 
\mbox{\rm in}\ \Omega\times
(0,+\infty)\label{Wo.1}\\
& &u (x,t) =0\quad \mbox{\rm on}\quad\partial\Omega\times
(0,+\infty)\label{Wo.2}\\
& &u(x,0)=u_0(x)\quad \mbox{\rm and}\quad u_t(x,0)=u_1(x)\quad \hbox{\rm
in}\quad\Omega\label{Wo.3}\\
& &u_t(x,t-\tau )=f(x,t),\quad \hbox{\rm
in}\quad\omega_2\times (0,\tau),\label{Wo.4}
\end{eqnarray}
with
initial
data $(u_0, u_1, f)\in H^1_0(\Omega)\times L^2(\Omega)\times L^2(\omega_2\times (0,\tau )),$  
and $a, k$ real constants, $a>\vert k\vert.$
The constant $\beta >0$ satisfies a suitable restriction to be specified below.

This problem enters into our previous framework, if we take
$H=L^2(\Omega)$ and the operator  $A_1$ defined by
$$A_1:{\mathcal D}(A_1)\rightarrow H\,:\,  u\rightarrow -\Delta u,$$
where ${\mathcal D}(A_1)=H^1_0(\Omega)\cap
H^2(\Omega).$  

The operator $A_1$ is a self--adjoint and positive operator with a compact inverse in $H$
and is such that $V={\mathcal D}(A^{1/2})=H^1_0(\Omega).$
We then define $W_i=L^2(\omega_i)$ and the operators $B, C$ as
$$
B:W_2\rightarrow H: \quad v\rightarrow \sqrt{k}\tilde v \chi_{\omega_2},
$$
$$
C:W_1\rightarrow H: \quad v\rightarrow \sqrt{a}\tilde v \chi_{\omega_1},
$$
where $\tilde v\in L^2(\Omega)$ is the extension of $v$ by zero outside $\omega_i.$ 
It is easy to verify that 
$$B^*(\varphi )=\sqrt {k} \varphi_{\vert_{\omega_2}}\quad\mbox{\rm for}\ \varphi\in H,$$
and thus 
$BB^*(\varphi)=k\varphi \chi_{\omega_2},$ for $\varphi\in H.$
Analogously,
$$C^*(\varphi )=\sqrt {a} \varphi_{\vert_{\omega_1}}\quad\mbox{\rm for}\ \varphi\in H,$$
and 
$CC^*(\varphi)=a\varphi \chi_{\omega_1},$ for $\varphi\in H.$
Moreover, since $\omega_2\subseteq \omega_1$ and $a>\vert k\vert $ the inequality $(\ref{contrast})$
holds.
Next, consider the functional 
$$G(u):=\frac 1 {\beta +2}\int_\Omega \vert u(x)\vert^{\beta+2} dx, \quad u\in H^1_0(\Omega ),$$ which, for $0<\beta \le \frac 4 {n-2},$ is well--defined by 
Sobolev's embedding theorem. 
Note that $G$ is G\^{a}teaux differentiable at any $u\in H^1_0(\Omega )$ and
its  G\^{a}teaux derivative is given by

$$DG (u)(v)=\int_{\Omega}\vert u(x)\vert^{\beta}u(x)v(x) dx,\quad v\in H^1_0(\Omega).$$

As proved in \cite{ACS}, 
if we assume that
$0< \beta <\frac 2 {n-2},$
then
$G$ satisfies the previous assumptions a), b), c).  
Therefore  problem (\ref{Wo.1})--(\ref{Wo.4}) enters in the abstract framework
(\ref{second}) and so the previous stability result holds for small initial data 
if the delay parameter $\vert k\vert$ is sufficiently small.
}
\end{Example}

\section{The case $\caB$ unbounded}

\label{ubdpbform}\hspace{5mm}

\setcounter{equation}{0}

In this case we need more assumptions on $\caB$, indeed we assume that
$$
\caB=\caC \caC^*,
$$
with
$\caC^*\in \caL(D(\caA), \caU)$ and hence $\caC\in \caL(\caU,  \caH_{-1})$,
where $\caU$ is    a complex
Hilbert space (which is identified with its dual space) and $\caH_{-1}=D(\caA^*)'$ is the dual space of $D(\caA^*)$ with respect to the pivot space $\caH$ (see \cite[section 2.10]{TucsnakWeiss_09}).
In such a setting,  for all $t\geq 0$ we can define $\Phi_t\in \caL(L^2(0,\tau; \caU), \caH_{-1})$ by
$$
\Phi_t v=\int_0^t S_{t-\sigma} \caC v(\sigma)\,d\sigma,
$$
for all $v\in L^2(0,\tau; \caU)$.

We further need the following assumptions that are satisfied by different examples (see below).
\\
{\bf (H1)} For all $v\in L^2(0,\tau; \caU)$, one has $\Phi_t v \in C([0,\tau], \caH)$ and there exists  $C_1>0$  such that
\be\label{h1}
\|\Phi_\tau v\|_\caH\leq C_1 \|v\|_{L^2(0,\tau; \caU)}.
\ee
By Remark  4.2.3 of \cite{TucsnakWeiss_09} we see that this property {\bf (H1)} implies that   $\caC$ is an admissible control operator for the semigroup generated by $\caA$ in the sense of Definition 4.2.1 of \cite{TucsnakWeiss_09}. Note further that
the estimate (4.2.5) of \cite{TucsnakWeiss_09} implies that
\be\label{h1'}
\|\Phi_t v\|_\caH\leq C_1 \|v\|_{L^2(0,\tau; \caU)}, \forall t\in [0,\tau].
\ee
{\bf (H2)} For all $v\in L^2(0,\tau; \caU)$, one has $\caC^*\Phi_t v \in L^2(0,\tau;  \caU)$ and there exists  $C_3>0$  such that
\be\label{h2}
\|\caC^*\Phi_t v\|_{L^2(0,\tau;  \caU)}\leq C_3 \|v\|_{L^2(0,\tau; \caU)}.
\ee
{\bf (H3)} For all $z_0\in \caH$, $\caC^* S(\cdot )z_0\in L^2(0,\tau; \caU)$
and
\be\label{h3}
\|\caC^* S(\cdot )z_0\|_{L^2(0,\tau;  \caU)}\leq C_2 \|z_0\|_{\caH}.
\ee

Note that this last condition directly implies that for all $\ell\in \N$,

\be\label{h3bis}
\|\caC^* S(\cdot )z_0\|_{L^2(\ell \tau,(\ell+1) \tau;  \caU)}\leq C_2 M e^{-\ell \tau \omega}\|z_0\|_{\caH}.
\ee
Indeed for any $t\in (\ell \tau,(\ell+1) \tau),$  we can write
$S(t)z_0=S(t-\ell \tau)S(\ell \tau)z_0,
$
and therefore by (\ref{h3})
\[
\|\caC^* S(\cdot )z_0\|_{L^2(\ell \tau,(\ell+1) \tau;  \caU)}
\leq C_2 \|S(\ell \tau) z_0\|_{\caH},
\]
which leads to (\ref{h3bis}) owing to our assumption
(\ref{assumpexpdecay}).

We are first able to prove the next well--posedness result.

\begin{Proposition}\label{ubdWellP}
Under the previous assumptions on $\caC$, then
for any initial datum $U_0\in {\mathcal H}$ and $f\in L^2(0,\tau; \caU)$, there exists a unique (mild) solution
$U\in C([0,+\infty), {\mathcal H})$  of problem
 \begin{eqnarray}\label{ubdabstract}
\left\{
\begin{array}{ll}
& U_{t}(t) = \caA U(t)+k \caC\caC^*  U(t-\tau)\quad \mbox{\rm in }
(0,+\infty)\\
& U(0)=U_0,\quad \caC^* U(t-\tau)=f(t), \quad \forall t\in (0,\tau).
\end{array}
\right.
\end{eqnarray}
\end{Proposition}

\noindent{\bf Proof.}
We use an iterative argument.
Namely in the interval $(0,\tau)$, problem (\ref{ubdabstract}) can be seen as an inhomogeneous evolution problem
\begin{eqnarray}\label{ubdabstract2}
\left\{
\begin{array}{ll}
& U_{t}(t) = \caA U(t)+k  \caC f(t)\quad \mbox{\rm in }
(0,\tau)\\
& U(0)=U_0.
\end{array}
\right.
\end{eqnarray}
Hence by  the hypothesis {\bf (H1)}, this problem has a unique
solution
$U\in  C([0,\tau], {\mathcal H})$ given by
\be\label{serge:19/02:2}
U(t)=S(t)U_0+k \int_0^t S(t-s) \caC f(s) \,ds.
\ee
This yields $U$ on $(0,\tau)$
and therefore on $(\tau,2\tau)$, problem (\ref{ubdabstract}) can be seen as an inhomogeneous evolution problem

\begin{eqnarray}\label{ubdabstract3}
\left\{
\begin{array}{ll}
& U_{t}(t) = \caA U(t)+k \caC v(t)\quad \mbox{\rm in }
(\tau, 2\tau)\\
& U(\tau)=U(\tau-),
\end{array}
\right.
\end{eqnarray}
where $v(t)= \caC^*  U(t-\tau)=\caC^*S(t-\tau)U_0+k \caC^*\Phi_{t-\tau} f$. But owing to the hypotheses
({\bf (H2)} and {\bf (H3)},
$v$ belongs to $L^2(\tau,2\tau; \caU)$.  Hence by  the hypothesis {\bf (H1)}, this problem has a unique solution
$U\in C([\tau,2\tau], {\mathcal H})$ given by
\be\label{serge:19/02:3}
U(t)=S(t-\tau)U(\tau-)+k\int_\tau^t S(t-s) \caC v(s)\,ds, \forall t\in [\tau,2\tau].
\ee
By iteration, we obtain a global solution. \qed

Similarly we will prove the following exponential stability result.

\begin{Theorem}\label{ubdstab2}
Let the assumptions {\bf (H1)} to  {\bf (H3)} be satisfied.
Set
\[
k_0:=\frac{e^{\tau\omega}-1}{M'^2 C_1 C_4  e^{2\omega \tau}},
\]
where $M'=\max\{M,1\}$, $C_4=\max\{C_2,\frac{C_3}{M' C_1}\}$.
Then  for any $k$ satisfying
\begin{equation}\label{ubdcondkpetit}
\vert k\vert < k_0,
\ee
 there exist $\omega' >0$ and $M''>0$ such that
the solution $U\in C([0,+\infty), {\mathcal H})$  of problem $(\ref{ubdabstract})$
with $U_0\in {\mathcal H}$ and $f\in L^2(0,\tau; \caU)$ satisfies
\begin{equation}\label{ubdexponentiald}
\|U(t)\|_\caH\le M'' e^{-\omega' t} (\|U_0\|_\caH+|k|   M' C_1 e^{2\omega \tau} \| f\|_{L^2(0,\tau, \caU)}),
\quad \forall t>0.
\end{equation}
From its definition the constant $k_0$ depends only on $M, \omega, \tau$ and the constants appearing in the assumptions {\bf (H1)} to  {\bf (H3)}.
\end{Theorem}
\noindent{\bf Proof.}
We use again an iterative argument and  the estimates (\ref{h1}) to (\ref{h3}).

First on $(0,\tau)$ using (\ref{serge:19/02:2}), the assumptions (\ref{assumpexpdecay}) and  (\ref{h1}), we see that
\[
\|U(t)\|_\caH\leq M e^{-\omega t} \|U_0\|_\caH+|k| C_1\| f\|_{L^2(0,\tau, \caU)}, \forall t\in (0,\tau),
\]
that directly  leads to
\be\label{serge:19/02:4}
\|U(t)\|_\caH\leq   e^{-\omega t} (M\|U_0\|_\caH+|k| C_1\| f\|_{L^2(0,\tau, \caU)} e^{\omega \tau}), \forall t\in (0,\tau).
\ee
Now coming back to (\ref{serge:19/02:2}) and using (\ref{h2}) and (\ref{h3}), we get
\be\label{serge:19/02:5}
\|\caC^* U\|_{L^2(0,\tau, \caU)}\leq C_2\|U_0\|_\caH+|k|C_3 \| f\|_{L^2(0,\tau, \caU)}.
\ee

Let us now prove by iteration that
for all $\ell\in \N$, we have
\be\label{serge:19/02:4l}
\|U(t)\|_\caH\leq   K_1(\ell) e^{-\omega t}, \forall t\in (\ell\tau,(\ell+1)\tau),
\ee
as well as
\be\label{serge:19/02:5l}
\|\caC^* U\|_{L^2(\ell\tau,(\ell+1)\tau; \caU)}\leq K_2(\ell) e^{- \ell \tau \omega},
\ee
where
\be\label{serge:19/02:k1}
K_1(\ell)\leq M'(\|U_0\|_\caH+\delta \alpha)(1+\delta C_4M')^{\ell},
\ee
and
\be\label{serge:19/02:k2}
K_2(\ell)\leq C_4M'(\|U_0\|_\caH+\delta \alpha)(1+\delta C_4M')^{\ell},
\ee
with   $\delta=|k| C_1 M' e^{2\omega \tau}$ and $\alpha=\| f\|_{L^2(0,\tau, \caU)}$.

Note that (\ref{serge:19/02:4l}) and (\ref{serge:19/02:5l}) hold for $\ell=0$ due to
(\ref{serge:19/02:4}) and (\ref{serge:19/02:5}) since simple calculations yield
\beqs
M\|U_0\|_\caH+|k| C_1\| f\|_{L^2(0,\tau, \caU)} e^{\omega \tau}\leq K_1(0),\\
C_2\|U_0\|_\caH+|k|C_3 \| f\|_{L^2(0,\tau, \caU)}\leq K_2(0).
\eeqs

Let us now prove that if (\ref{serge:19/02:4l})--(\ref{serge:19/02:k2}) hold up to  $\ell$
then they hold for $\ell+1$.

Indeed for $t\in ((\ell+1)\tau,(\ell+2)\tau)$, we have
\beqs
U(t)&=&S(t)U_0+
k\sum_{j=1}^{\ell}\int_{j\tau}^{(j+1)\tau}S(t-s)\caC\caC^*U(s-\tau)\,ds
\\
&+&
k\int_{(\ell+1)\tau}^tS(t-s)\caC\caC^*U(s-\tau)\,ds
\nonumber
\\
&+&
k\int_0^\tau S(t-s)\caC f(s)\,ds.
\nonumber
\eeqs
This identity can be equivalently written
\beq\label{serge19/02:10}
U(t)&=&S(t)U_0+
k \sum_{j=1}^{\ell} S(t-(j+1)\tau)\Phi_\tau \caC^*U((j-1)\tau+\cdot)
\\
&+&
k\Phi_{t- (\ell+1)\tau} \caC^*U(\ell\tau+\cdot)
\nonumber
+
k S(t-\tau)\Phi_\tau   f.
\nonumber
\eeq
Hence by our assumptions (\ref{assumpexpdecay}) and  (\ref{h1}), we deduce that
\beqs
\|U(t)\|_\caH&\leq &M e^{-\omega t} \|U_0\|_\caH+|k|
 M C_1 \sum_{j=1}^{\ell} e^{-(t-(j+1)\tau)\omega} \|\caC^*U((j-1)\tau+\cdot)\|_{L^2(0,\tau,\caU)}
\\
&+&
|k|C_1\|\caC^*U(\ell\tau+\cdot)\|_{L^2(0,\tau,\caU)}
\\
&+&
|k|M e^{-(t-\tau)\omega} C_1 \alpha.
\eeqs
Hence by our iterative assumption, the estimate (\ref{serge:19/02:5l}) for all $j\leq \ell$
yields
\beqs
\|U(t)\|_\caH&\leq &M' e^{-\omega t} \Big (\|U_0\|_\caH+|k|
  C_1 e^{2\tau\omega}\sum_{j=1}^{\ell+1}   K_2(j-1)
\\
&+&
|k|  e^{\tau\omega} C_1 \alpha\Big).
\eeqs
By setting $K_2(-1)=\alpha$, we have found that
\beqs
\|U(t)\|_\caH&\leq &M' e^{-\omega t} \Big(\|U_0\|_\caH+|k|
  C_1  e^{2\tau\omega}\sum_{j=0}^{\ell+1}   K_2(j-1)\Big).
\eeqs
This proves (\ref{serge:19/02:4l}) for $\ell+1$
with
\be\label{serge:19/02:k1iter}
K_1(\ell+1)= M'\Big(\|U_0\|_\caH+|k|
  C_1  e^{2\tau\omega}\sum_{j=0}^{\ell+1}   K_2(j-1)\Big).
\ee

Now we come back to (\ref{serge19/02:10}) and applying $\caC^*$ to this identity (meaningful due to our assumptions {\bf (H2)} and {\bf (H3)}), we get

\beqs
\|\caC^*U\|_{L^2((\ell+1)\tau,(\ell+2)\tau, \caU)}&\leq&C_2 M e^{-(\ell+1)\tau)\omega} \|U_0\|_\caH\\
&+&
|k|C_2C_1M\sum_{j=1}^{\ell} e^{-(\ell-j)\tau\omega}\|\caC^*U((j-1)\tau+\cdot)\|_{L^2(0,\tau, \caU)}
\\
&+&
|k|C_3\|\caC^*U \|_{L^2(\ell\tau,(\tau+1)\tau, \caU)}
\nonumber
+
|k|C_2C_1 e^{-\ell \tau\omega}\alpha.
\nonumber
\eeqs
As our iterative assumption means that (\ref{serge:19/02:5l}) holds for all $j\leq \ell$, we get
\beqs
\|\caC^*U\|_{L^2((\ell+1)\tau,(\ell+2)\tau, \caU)}&\leq&C_2 M e^{-(\ell+1)\tau\omega} \|U_0\|_\caH+
|k|C_2C_1M' e^{-\ell \tau\omega}e^{\tau\omega}\sum_{j=0}^{\ell}  K_2(j-1)
\\
&+&
|k|C_3e^{-\ell \tau\omega} K_2(\ell).
\eeqs
As $C_2 \leq  C_4$ and $C_3\leq C_1M'C_4$, we deduce that
\beqs
\|\caC^*U\|_{L^2((\ell+1)\tau,(\ell+2)\tau, \caU)}&\leq&C_4 M' e^{-(\ell+1)\tau\omega} \Big(\|U_0\|_\caH+
\delta \sum_{j=-1}^{\ell}  K_2(j)\Big).
\eeqs
This proves (\ref{serge:19/02:5l}) for $\ell+1$
with
\be\label{serge:19/02:k2iter}
K_2(\ell+1)= C_4M'\big(\|U_0\|_\caH+\delta \sum_{j=-1}^{\ell}  K_2(j)\big).
\ee

Let us now show that $K_2(\ell)$ given by (\ref{serge:19/02:k2iter}) satisfies (\ref{serge:19/02:k2}).
Indeed it holds for $\ell=0$ and then we again prove (\ref{serge:19/02:k2}) by induction.
If it holds up $\ell$ then by  (\ref{serge:19/02:k2iter}) we will have
\beqs
K_2(\ell+1)&\leq& C_4M'\Big(\|U_0\|_\caH+\delta \alpha +\delta \sum_{j=0}^{\ell}  C_4M'(\|U_0\|_\caH+\delta \alpha)(1+\delta C_4M')^j\Big)
\\
&\leq& C_4M'(\|U_0\|_\caH+\delta \alpha )\Big(1+\delta C_4M' \sum_{j=0}^{\ell}   (1+\delta C_4M')^j\Big)
\\
&\leq& C_4M'(\|U_0\|_\caH+\delta \alpha )\Big(1+\delta C_4M' \frac{(1+\delta C_4M')^{\ell+1}-1}{\delta C_4M'}\Big)
\\
&\leq& C_4M'(\|U_0\|_\caH+\delta \alpha ) (1+\delta C_4M')^{\ell+1}.
\eeqs
This proves (\ref{serge:19/02:k2}) for $\ell+1$.

Once (\ref{serge:19/02:k2}) holds for all $\ell$, we come back to (\ref{serge:19/02:k1iter}) and get
\beqs
K_1(\ell+1)&\leq& M'\big(\|U_0\|_\caH
\\
&+&|k|
  C_1  e^{2\tau\omega}(\alpha+\sum_{j=1}^{\ell}   C_4M'(\|U_0\|_\caH+\delta \alpha)(1+\delta C_4M')^{j})\Big)
\\
&\leq& M'\big(\|U_0\|_\caH
\\
&+&\delta
(\alpha+\delta^{-1}(\|U_0\|_\caH+\delta \alpha)((1+\delta C_4M')^{\ell+1}-1)
\Big)
\\
&\leq& M'\big(\|U_0\|_\caH
\\
&+&
(\alpha\delta+(\|U_0\|_\caH+\delta \alpha)((1+\delta C_4M')^{\ell+1}-1)
\Big)
\\
&\leq& M'\Big(\|U_0\|_\caH
+
(\|U_0\|_\caH+\delta \alpha)(1+\delta C_4M')^{\ell+1}-\|U_0\|_\caH\Big)
\\
&\leq& M' (\|U_0\|_\caH+\delta \alpha)(1+\delta C_4M')^{\ell+1}.
\eeqs
This proves (\ref{serge:19/02:k2}) for all $\ell+1$.

We end up the proof by combining (\ref{serge:19/02:4l}) and (\ref{serge:19/02:k1}) to get

\beqs
\|U(t)\|_\caH\leq   M'(\|U_0\|_\caH+\delta \alpha)(1+\delta C_4M')^{\ell}  e^{-\omega t}, \forall t\in (\ell\tau,(\ell+1)\tau).
\eeqs
Hence setting
\begin{equation}\label{ubddefsigma}
\sigma=\tau^{-1} \ln(1+\delta C_4M')=\tau^{-1} \ln(1+|k| C_4M'^2 C_1  e^{2\omega \tau}),
\end{equation}
we conclude as in Theorem \ref{stab2} that
\begin{equation}\label{ubdserge14/10:1}
\|U(t)\|_\caH\leq  M'(\|U_0\|_\caH+\delta \alpha)  e^{(\sigma-\omega) t}, \forall t>0.
\end{equation}

Therefore $\|U(t)\|_\caH$ will decay exponentially if
$\sigma-\omega$ is negative or equivalently if
\[
1+|k| C_4M'^2 C_1  e^{2\omega \tau}<e^{\tau\omega},
\]
which is nothing else than (\ref{ubdcondkpetit}).
\qed

\section{Examples in the case $\caB$ unbounded}

Most of our examples are second order evolution equations with damping.
Namely they are in the following  form. Let $H$ be a complex Hilbert space and let $A:{\mathcal D}(A)\rightarrow H$
be a positive self--adjoint  operator with a compact inverse in $H.$ Denote by $V:={\mathcal D}(A^{\frac 1 2})$ the domain of
$A^{\frac 1 2}.$ Moreover,
for $i=1,2,$
let $U_i$ be
     complex Hilbert spaces
with norm and inner product denoted respectively by $\Vert \cdot\Vert_{U_i}$
and $\langle\cdot ,\cdot\rangle_{U_i}$
and let $B_i :U_i\rightarrow V^\prime$          be
linear
operators. In this setting we consider the problem
 \begin{eqnarray}
& &u_{tt}(t) +A u (t)+B_1 B_1^* u_t(t) +k B_2 B_2^* u_t(t-\tau) =0\quad t>0\label{1.1}\\
& &u(0)=u_0\quad \mbox{\rm and}\quad u_t(0)=u_1\quad\label{1.2}
\end{eqnarray}
where the constant  $\tau >0$ is the time delay and $k$ is a real parameter.

We transform this problem into a first order system by using the standard reduction of order:
setting
$$
U=(u,  u_t)^\top,
$$
it satisfies formally
$$
U_t={\cal A} U-k\caC\caC^* U(t-\tau) ,\  U(0)=(u_0,  u_1)\in V\times H,
$$
where
$$
{\cal A}  (u,v)^\top=(v,   -Au- B_1 B_1^* v),
$$
with
$$
D({\cal A})=\{(u,v)^\top\in V\times V: Au+  B_1 B_1^* v\in H\},
$$
and
$$
\caC=\left(
\begin{array}{ll}
0&0\\
0&B_2
\end{array}
\right).
$$

In such a setting, we easily check that the adjoint
$\caA^*$ of $\caA$ is given by
$$
{\cal A}^* (u,v)^\top=(-v,   Au- B_1 B_1^* v),
$$
with
$$
D({\cal A}^*)=\{(u,v)^\top\in V\times V: Au-  B_1 B_1^* v\in H\}.
$$
In other words, if we introduce the unitary mapping
$$
O=\left(
\begin{array}{ll}
Id&0\\
0&-Id
\end{array}
\right),
$$
we see that
$$
{\cal A}^* (u,v)^\top={\cal A} O (u,v)^\top.
$$
Consequently the semigroup $(S^*(t))_{t\geq 0}$ generated by ${\cal A}^*$    will be given by
\[
S^*(t)=OS(t) O.
\]

To apply our stability results from section \ref{ubdpbform} to our system (\ref{1.1})--(\ref{1.2}) we need to check the assumptions {\bf (H1)} to {\bf (H3)} for the operators $\caA$ and $\caB=\caC\caC^*$.
But in this case, {\bf (H1)}  implies {\bf (H3)}  since by Remark 4.2.4 of \cite{TucsnakWeiss_09},
 {\bf (H1)}  implies that $\caC$ is an admissible control for the semigroup $S(t)$
and by Theorem 4.4.3 of \cite{TucsnakWeiss_09} this is equivalent to the fact that
$\caC^*$ is an admissible operator for the semigroup   $S^*$. As
$\caC^*S^*=O\caC^* SO$, we deduce that  {\bf (H3)} holds owing to Proposition  4.4.1 of \cite{TucsnakWeiss_09}.

\subsection{The wave equation with boundary feedbacks in 1d\label{ss1d}}

Our first application concerns the wave equation with boundary feedbacks in dimension 1.
More precisely let $\Omega=(0,1)\subset\RR$ be the unit interval.

Given  a  positive constant  $a$,
let us consider the initial boundary value  problem

 \begin{eqnarray}
& &u_{tt}(x,t) -  u_{xx} (x,t)=0\quad \mbox{\rm in}\quad\Omega\times
(0,+\infty)\label{W.1}\\
& &  u_x (1,t)=-  u_t(1,t)-au(1,t)\quad \mbox{\rm on}\quad
(0,+\infty)\label{W.2B}\\
& &  u_x (0,t)=k  u_t(0,t-\tau)\quad \mbox{\rm on}\quad
(0,+\infty)\label{W.2Bdelay}\\
& &u(x,0)=u_0(x)\quad \mbox{\rm and}\quad u_t(x,0)=u_1(x)\quad \hbox{\rm
in}\quad\Omega\label{W.3}
\end{eqnarray}
with
initial
data $(u_0, u_1)\in H^1(\Omega)\times L^2(\Omega)$.
This problem enters in the abstract framework (\ref{1.1})--(\ref{1.2}) if we take
$H=L^2(\Omega)$ with its standard norm and  $V=H^1(\Omega )$ with the norm
\[
\|u\|_V^2=\int_0^1|u_x(x)|^2\,dx+a |u(1)|^2, \forall u\in H^1(\Omega ).
\]
The operator  $A$ is defined by
$$A:{\mathcal D}(A)\rightarrow H\,:\,  u\rightarrow -\Delta u,$$
where
$${\mathcal D}(A):=\{u\in H^2(\Omega)\cap H^1(\Omega )
\,:\,   u_x(0)=0 \hbox{ and } u_x(1)+a u(1)=0\}.$$

We then define $U_1=U_2:=\RR$
and for $i=1$ or 2, the operator $B_i^*\in {\cal L}(V,U_i)$ as
$$B_1^*w =   w(1),\  B_2^*w =   w(0),\,\  \forall \,w \in
V.$$
Consequently
$B_i \in {\cal L}(U_i;V^\prime    )$ is characterized as follows: for any $v \in \RR$,
$$B_1v=  v\delta_1, B_2v=\  v\delta_0.$$
Finally we need to take $\caU=\{0\}\times U_2$.

Hence in such a situation it remains to check  the hypotheses {\bf (H1)} and {\bf (H2)}.

To check the assumption  {\bf (H2)},  as $\caD(0,\tau)$ is dense in $L^2(0,\tau)$,
it suffices to check it for   $v\in \caD(0,\tau)$. For such a $v$  consider
$u=\Phi_t v$, $0<t<\tau$, that is the (strong) solution
  of
 \begin{eqnarray}\label{pbphit}
\left\{
\begin{array}{llll}
& &u_{tt}(x,t) - u_{xx} (x,t)=0\quad \mbox{\rm in}\quad\Omega\times
(0,\tau), \\
& &  u_x (1,t)=- u_t(1,t)-au(1,t)\quad \mbox{\rm on}\quad
(0,\tau),\\
& &u_x (0,t)=- v(t)\quad \mbox{\rm on}\quad
(0,\tau), \\
& &u(x,0)=0\quad \mbox{\rm and}\quad u_t(x,0)=0\quad \hbox{\rm
in}\quad\Omega.
\end{array}
\right.
\end{eqnarray}
Consider an extension $\tilde v$ of $v$ by taking an odd extension of $v$ to $(\tau, 2\tau)$ and by taking $\tilde v=0$
outside $(0,2\tau)$. Then such an extension satisfies
$$
\int_\RR \tilde v(t)\,dt=0,
$$
and
\be\label{serge19/02:1}
\int_\RR |\tilde v(t)|^2\, dt\leq 2\int_0^\tau |v(t)|^2\, dt.
\ee
Then we can consider the solution $w$ of
\begin{eqnarray*}
& &w_{tt}(x,t) -  w_{xx} (x,t)=0\quad \mbox{\rm in}\quad\Omega\times
(0,\infty), \\
& &  w_x (1,t)=- w_t(1,t)-aw(1,t)\quad \mbox{\rm on}\quad
(0,+\infty),\\
& &w_x (0,t)=- \tilde v(t)\quad \mbox{\rm on}\quad
(0,+\infty), \\
& &w(x,0)=0\quad \mbox{\rm and}\quad w_t(x,0)=0\quad \hbox{\rm
in}\quad\Omega.
\end{eqnarray*}
But since the corresponding operator $\caA$ generates a strongly continuous semigroup,
this solution $w$ coincides with $u$ in $(0,\tau)$. Furthermore we can extend $w$ by zero in $(0,1)\times (-\infty,0)$ that  then satisfies
\begin{eqnarray*}
& &w_{tt}(x,t) - w_{xx} (x,t)=0\quad \mbox{\rm in}\quad\Omega\times
\RR, \\
& &  w_x (1,t)=- w_t(1,t)-aw(1,t)\quad \mbox{\rm on}\quad
\RR,\\
& &w_x (0,t)=- \tilde v(t)\quad \mbox{\rm on}\quad
\RR.
\end{eqnarray*}
Taking Fourier transform in time, we deduce that for all $\xi\in \RR$, $\hat w(\cdot,\xi)$ satisfies
\begin{eqnarray*}
& &\xi^2 \hat w + \hat w_{xx}  =0\quad \mbox{\rm in}\quad\Omega, \\
& &  \hat w_x (1)=-(a+i\xi) \hat w (1),\\
& &\hat w_x (0)=- \widehat{\tilde v}(\xi).
\end{eqnarray*}

Hence     easy calculations show that%ici
\[
\hat w(x,\xi)= \frac{\widehat{\tilde v}(\xi)}{i\xi} (e^{-i\xi x}+c(\xi) \cos(\xi x)), \ \forall x\in (0,1),
\]
with
\[
c(\xi)=\frac{a  e^{-i\xi}}{\xi\sin \xi-(a+i\xi)\cos\xi}.
\]

This identity implies that
\[
i\xi \hat w(0,\xi)=  \widehat{\tilde v}(\xi) (1+c(\xi))
\]
and since one can show that there exists a positive constant $C$ depending on $a$ such that
\[
|c(\xi)|\leq C, \ \forall \xi\in \RR,
\]
we deduce that
\[
|i\xi \hat w(0,\xi)|\leq (1+C)  |\widehat{\tilde v}(\xi)|, \forall \xi\in \RR.
\]
By Parseval's identity we find that
$$\int_\RR |w_t(0, t)|^2\,dt\leq (1+C)
\int_\RR |\tilde v (t)|^2\,dt.
$$
Recalling that $w$ coincides with $u$ in $(0,\tau)$ and using the estimate (\ref{serge19/02:1}),
we have proved that
\be\label{serge16/04:1}
\int_0^\tau  |u_t(0, t)|^2\,dt\leq 2 (1+C)
\int_0^\tau |v (t)|^2\,dt.
\ee
This implies that
{\bf (H2)} holds reminding that
$$
\caC^* (u,u_t)=B_2^* u_t=u_t(0,\cdot).
$$

As before it suffices to check the  assumption  {\bf (H1)} for $v\in \caD(0,\tau)$. For such a $v$  consider
   the (strong) solution $u=\Phi_t v$, $0<t<\tau$
  of (\ref{pbphit}).
Then we consider its energy
\[
\caE(t)=\frac12(\int_0^1(|u_t|^2+|u_x|^2)\,dx+a |u(1,t)|^2).
\]
Differentiating and integrating by parts we have
\[
\caE'(t)=-|u_t(1,t)|^2+v(t) u_t(0,t).
\]
Integrating this identity between $0$ and $t\in (0,\tau]$ and using Cauchy-Schwarz's inequality we find that
\[
\caE(t)\leq \|v\|_{L^2(0,\tau)} \|u_t(0,t)\|_{L^2(0,\tau)}.
\]
Hence using the estimate (\ref{serge16/04:1}) we arrive at (\ref{h1}).

The continuous property is proved similarly by integrating between $t\in (0,\tau]$ and $t'\in (0,\tau]$.

In conclusion, as the system (\ref{W.1})--(\ref{W.3}) with $k=0$
is exponentially stable and the assumptions {\bf (H1)} to {\bf (H3)} hold,
 system (\ref{W.1})--(\ref{W.3})
remains exponentially stable if   $k$ is small enough.

\br
{\rm
Our approach cannot be used for the wave equation in $\RR^d$, with $d\geq 2$
since according to the results from \cite{Tataru:98} (see for instance Theorem 3 in \cite{Tataru:98}
and the comments before), the assumption  {\bf (H2)} is wrong once  $d\geq 2$.
}
\er

\subsection{The wave equation with boundary and internal unbounded feedbacks in 1d}

Here we want to consider the following problem: For a fixed $a\in (0,1)$ consider the solution of
\begin{eqnarray}
& &u_{tt}(x,t) - u_{xx} (x,t)=0\quad \mbox{\rm in}\quad(0,a)\cup (a,1)\times
(0,+\infty)\label{W.1-1d}\\
& &u (0,t) =0\quad \mbox{\rm on}\quad
(0,+\infty)\label{W.2-1d}\\
& &  u_x (1,t)=- u_t(1,t)\quad \mbox{\rm on}\quad
(0,+\infty)\label{W.2B-1d}\\
& &[u](a)=0,\quad [u_x](a)=k u_t(a,t-\tau)\quad \mbox{\rm on}\quad
(0,+\infty)\label{W.2Bdelay-1d}\\
& &u(x,0)=u_0(x)\quad \mbox{\rm and}\quad u_t(x,0)=u_1(x)\quad \hbox{\rm
in}\quad(0,1)\label{W.3-1d}
\end{eqnarray}
with
initial
data $(u_0, u_1)\in \{w\in H^1(0,1): w(0)=0\}\times L^2(0,1),$  and $[u](a)$ means the jump of $u$ at the point $a$, i.e., $[u](a)=u(a+)-u(a-)$. This problem corresponds to the case where a standard dissipative law (cfr. \cite{Komornikbook}) is acting at 1, while a dissipation with delay appears at the interior point $a$.

As in subsection \ref{ss1d}, we only need to check the assumption {\bf (H2)} (since as before one can show that {\bf (H2)} implies {\bf (H1)}),
that is proved exactly as before by using an extension method and Fourier transform in time to get the system
\begin{eqnarray*}
& &\xi^2 \hat w - \hat w_{xx}  =0\quad \mbox{\rm in}\quad (0,a)\cup (a,1), \\
& &\hat w (0)=0,\\
& &  \hat w_x (1)=-i\xi \hat w (1),\\
& &[\hat w](a)=0,\quad [\hat w_x](a)= \widehat{\tilde v}(\xi).
\end{eqnarray*}
Again simple calculations yield for $\xi\ne 0$
\beqs
\hat w(x,\xi)=c_1\sin (\xi x)\quad \hbox{ in } (0,a),
\\
\hat w(x,\xi)=c_2\cos (\xi x)+c_3\sin (\xi x)\quad \hbox{ in } (a,1),
\eeqs
with $|c_1|=|\xi|^{-1}$ and $|c_2|=|c_3|\leq |\xi|^{-1}$.
This directly implies
\[
|\xi \hat w (a,\xi)|\leq 1,
\]
and leads to the conclusion because here
$B_2^*v=v(a)$.

In conclusion, for $k$ small enough, system (\ref{W.1-1d})--(\ref{W.3-1d}) is exponentially stable
since it is for $k=0$.

\subsection{The wave equation with a bounded internal  feedback and a boundary   unbounded feedback in 1d}

Arguing as before we can consider the following problem
\begin{eqnarray}
& &u_{tt}(x,t) - u_{xx} (x,t)+\alpha u_t(x,t)=0\quad \mbox{\rm in}\quad(0,1)\times
(0,+\infty)\label{W.1-1dint}\\
& &u (0,t)  =0\quad \mbox{\rm on}\quad
(0,+\infty)\label{W.2-1dint}\\
& &  u_x (1,t)=k u_t(1,t-\tau)\quad \mbox{\rm on}\quad
(0,+\infty)\label{W.2B-1dint}\\
& &u(x,0)=u_0(x)\quad \mbox{\rm and}\quad u_t(x,0)=u_1(x)\quad \hbox{\rm
in}\quad(0,1)\label{W.3-1dint}
\end{eqnarray}
with $\alpha>0$ and
initial
data $(u_0, u_1)\in \{w\in H^1(0,1): w(0)=0\}\times L^2(0,1).$ This problem corresponds to the case where a standard dissipative law is acting on the whole domain, while a dissipation with delay appear at the boundary point $1$.
As this system   with $k=0$
is exponentially stable and the assumptions {\bf (H1)} to {\bf (H3)} are valid,
 system (\ref{W.1-1dint})--(\ref{W.3-1dint})
remains exponentially stable if   $k$ is small enough (cfr. \cite{DLP}).

 \end{document}